\documentclass[submission]{FPSAC2021}

\newtheorem{thm}{Theorem}

\newtheorem{definition}{Definition}
\newtheorem{lem}{Lemma}

\usepackage{lipsum}
\usepackage{enumitem}

\usepackage{tikz-cd}
\usetikzlibrary{arrows}
\usetikzlibrary{decorations.pathmorphing}

\newcommand{\N}{\mathbb{N}}
\newcommand{\ka}{\kappa}

\title[LLT cumulants of unicellular diagrams]{LLT cumulants of unicellular Young diagrams, parking functions and Schur positivity}

\author[M. Kowalski]{Maciej Kowalski\thanks{\href{mailto:mkowalski@impan.pl}{mkowalski@impan.pl} MK is supported by {\it Narodowe Centrum Nauki}, grant UMO-2017/26/D/ST1/00186. }\addressmark{1}}

\address{\addressmark{1}Institute of Mathematics, Polish Academy of Sciences, Warszawa, Poland}

\received{\today}

\abstract{We give a combinatorial formula for LLT cumulants of unicellular unilevelled shapes in terms of parking functions and Cayley trees. Our formula implies previously conjectured Schur positivity.}


\usepackage[backend=bibtex,doi=false,isbn=false,url=false]{biblatex}
\DeclareMathOperator{\LLT}{LLT}
\DeclareMathOperator{\Inv}{Inv}
\DeclareMathOperator{\inv}{inv}
\DeclareMathOperator{\Part}{Part}
\DeclareMathOperator{\SSYT}{SSYT}

\begin{document}

\maketitle

\section{Introduction}

\textit{LLT polynomials} were introduced by Lascoux, Leclerc, and Thibon \cite{LLTDef} in the context of quantum groups. Their importance in the theory of symmetric functions quickly became apparent mostly due to the rich yet hidden combinatorics and a broad spectrum of relations --- for instance, with the celebrated Macdonald polynomials, Kazhdan--Lusztig theory, algebraic geometry, and knots theory.

\subsection{Connections to $e$-positivity}

In the following, we will mostly study the special case of \emph{unicellular LLT polynomials}. These are particularly interesting in the context of the $e$-positivity phenomenon, which recently gained a lot of attention. One of the most interesting open problems in this field is a conjecture of Stanley and Stembridge \cite{StanleyStembridge} on $e$-positivity of chromatic symmetric functions of $(3+1)$-free posets. After a breakthrough of Guay-Paquet \cite{GuayPaquetChrUnicell}, who reduced the Stanley--Stembridge conjecture to the case of unit-interval graphs, a lot of research appeared in this and related context. Among others, Alexandersson \cite{Alexandersson2020} conjectured an explicit $e$-positive formula for the unicellular LLT polynomials after the shift $q \longmapsto q+1$, which was recently proven by Alexandersson and Sulzgruber \cite{AlexSulz}\footnote{In fact, the $e$-positivity phenomenon holds true for a more general class of vertical-strip LLT polynomials, but this follows easily from the $e$-positivity of unicellular LLT polynomials}.

The relation between the unicellular LLT polynomials and unit interval graphs is given by the following observation: unicellular LLT polynomials are indexed by Dyck paths, and those are in a natural one-to-one correspondence with unit interval graphs $D \leftrightarrow G_D$ (see \cref{sec:LLTGraphs}). Therefore, after the shift $q \longmapsto q+1$ we have the following formula
\begin{equation} \label{def:lltdyck}
\LLT_{D}(q+1) = \sum_{H \subseteq G_{D}}\sum_{f:V(H)\to \N_+} q^{e(H)}\prod_{v \in V(H)}x_{f(v)},
\end{equation}
where we sum over all subgraphs $H$ and all the \emph{valid colorings} $f$ of $H$ --- those where $f(v) \geq f(w)$ whenever we have a directed edge $v \to w$ in $H$ (here $e(H)$ denotes the number of edges in $H$).

This decomposition suggests that the $e$-positivity phenomenon is not necessarily related to the abstract framework of LLT polynomials, but rather to this particular interpretation as generating functions of certain graphs. In particular, motivated by this observation, we introduce the following function, which can be interpreted as the ``building block'' of unicellular LLT polynomials after the shift $q \longmapsto q+1$ and which is the main object in this paper:
\begin{equation} \label{def:graphcum}
\ka_{D}(q) = \sum_{H \subseteq G_{D}}\sum_{f:V(H)\to \N_+}
(q-1)^{e(H)-v(H)+1}\prod_{v \in V(H)}x_{f(v)},
\end{equation}
where we sum over all {\bf connected} subgraphs $H$ of a graph $G_D$, and all the valid colorings $f$ of $H$ (here $e(H)$ and $v(H)$ denote the number of edges and vertices in $H$ respectively). We call this function the \emph{unicellular LLT cumulant}. In the joint paper with Dołęga~\cite{DolegaKowalski} we introduce LLT cumulants in the general case by using the formula for cumulants from the probability theory and we prove their combinatorial formula generalizing \eqref{def:graphcum} by replacing unit interval graphs by more general graphs. In particular, we explain that the $e$-positivity phenomenon indeed holds true for the vertical-shape LLT cumulants with the same combinatorial interpretation as in~\cite{AlexSulz}, which refines Alexandersson--Sulzgruber result and gives a combinatorial interpretation of the shift $q \longmapsto q+1$.

\subsection{Connections to Macdonald polynomials}

Although LLT cumulants appear naturally in the aforementioned context of the \linebreak $e$-positivitiy, our main motivation for introducing them comes from the Schur-positivity conjecture for Macdonald cumulants. The latter were introduced by Dołęga \cite{DolegaMacCum,Dolega2017} and, roughly speaking, can be interpreted as the higher-order approximation of the product of Macdonald polynomials as $q\to 1$. The trivial case --- a cumulant of order one --- corresponds to Macdonald polynomials, for which the Schur-positivity is the celebrated result of Haiman~\cite{HaimanMacdonaldpositivity}. Therefore the conjecture of Dołęga seems to be a challenging open problem, whose solution might shed new light on the connection between symmetric functions, representation theory, and algebraic geometry.

One of the beautiful results of Haiman, Haglund, and Loehr~\cite{HaimanLLTaltdef} is a formula expressing Macdonald polynomials as a positive linear combination of LLT polynomials. Together with the Schur-positivity of LLT polynomials~\cite{GrojnowskiHaiman2007}, they give an alternative proof of the Schur positivity of Macdonald polynomials. Inspired by this result, we found a formula~\cite{DolegaKowalski} for the expansion of Macdonald cumulants as a positive linear combination of LLT cumulants, and we posed several conjectures on the expansion of LLT cumulants as linear combinations of LLT polynomials, which refine the Schur positivity conjecture of Dołęga for Macdonald cumulants.

In the special case of unicellular LLT cumulants, we conjectured the following combinatorial formula:
\begin{equation} \label{conj:Tdecgen}
\ka_{D}(q) = \sum\limits_T \LLT(\nu(T)),
\end{equation}
where the sum runs over some spanning trees of $G_{D} $ and $\nu(T)$ denotes a sequence of vertical-strip shapes corresponding to $T$ (see \cref{sec:cumunicell}).

Note that formula \eqref{conj:Tdecgen} is a priori very different from the aforementioned $e$-expansion which generalizes the result of Alexandersson and Sulzgruber --- it implies Schur-positivity of $\ka_{D}(q)$, which is not clear from the $e$-expansion after the shift $q \longmapsto q+1$.

Our main result gives the proof of \eqref{conj:Tdecgen} in the special case of the maximal Dyck path $D_n$ which corresponds to the complete graph $G_{D_n} = K_n$. In this case

\begin{thm} \label{thm:singcelldecomp}
	\begin{equation} \label{thm:formulawithtrees}
	\ka_{D_n} = \sum_{T\subset K_n} \LLT(\nu(T)) = \sum_{f\in PF_{n-1}} \LLT(\mu(f)),
	\end{equation}
	where the first sum runs over all Cayley trees on $[n]$ vertices, the second over all parking functions on $n-1$ cars, and $\nu(T)$ and $\mu(f)$ are sequences of vertical-strip shapes associated to such objects (see \cref{sec:cumunicell} and \cref{sec:schroderpaths}).
\end{thm}

As an immediate corollary of \cref{thm:singcelldecomp}, we get the proof of the conjecture of Dołęga~\cite{DolegaMacCum} on Schur-positivity of Macdonald cumulants in the
special case $\lambda^1= \dots = \lambda^n = (1)$.

We would like to mention the similarity between \eqref{thm:formulawithtrees} and the celebrated formula for the character of the diagonal coinvariant algebra \cite{CarlssonMellitShuffle}. Both formulas are given as linear combinations of LLT polynomials associated to parking functions, but the LLT polynomials corresponding to the fixed parking function differs in the two. It would be interesting to find a representation-theoretical interpretation of LLT cumulants, which encourages further investigation.

The main idea for the proof is to switch from the standard framework, where the combinatorial objects of focus are Young diagrams, Dyck, or Schr\"oder paths, to a more general one: that of graphs with colored edges. We are going to study the generating functions of these graphs, and the price we pay is that these are no longer symmetric functions. Nevertheless, we prove they generate the algebra of quasi-symmetric functions, where we can classify all the relations by a variant of the result of F\'eray~\cite{FerayIncExc}. Using these, we are able to transform the generating function corresponding to the LHS of \eqref{thm:formulawithtrees} to generating functions corresponding to the RHS of \eqref{thm:formulawithtrees}. We believe that due to the flexibility of our framework, it might be of independent interest.

\section{LLT graphs} \label{sec:LLTGraphs}

Instead of the original definition of LLT polynomials \cite{LLTDef}, we use that of Haglund, Haiman, and Loehr \cite{HaimanLLTaltdef}.

Let $\pmb{\lambda}/\pmb{\mu}=(\lambda_1/\mu_1,\dots,\lambda_n/\mu_n)$ be a sequence of skew shapes. We say that a cell $\alpha=(i,j)$ of $\lambda_k/\mu_k$ has \emph{content} $c(\alpha) := i-j$ and \emph{shifted content} $\tilde{c}(\alpha) := nc(\alpha) + k$.

Let $\mathcal{T}\in\SSYT(\pmb{\lambda}/\pmb{\mu})$, where $\SSYT(\pmb{\lambda}/\pmb{\mu})$ denotes the set of semistandard Young tableaux on $\pmb{\lambda}/\pmb{\mu}$. The set of \emph{inversions} in $\mathcal{T}$ is defined to be $\Inv(\mathcal{T}) := \{(\alpha,\beta)\in\pmb{\lambda}/\pmb{\mu} \mid 0 < \tilde{c}(\beta)-\tilde{c}(\alpha) < n \text{ and } \mathcal{T} (\alpha)>\mathcal{T} (\beta)\}$. We denote $\inv(\mathcal{T}):=|\Inv(\mathcal{T})|$.

The \textit{LLT polynomial} of $\pmb{\lambda}/\pmb{\mu}$ is the generating function
\begin{equation}
\LLT(\pmb{\lambda}/\pmb{\mu}) := \sum_{\mathcal{T}\in\SSYT(\pmb{\lambda}/\pmb{\mu})} q^{\inv(\mathcal{T})}x^\mathcal{T}.
\end{equation}
If $\pmb{\lambda}/\pmb{\mu}$ is a sequence of unicellular shapes, then $\LLT(\pmb{\lambda}/\pmb{\mu})$ is called a \emph{unicelullar LLT polynomial}.

Following \cite{DolegaKowalski}, we define the unicellular \textit{LLT cumulant} to be the expression
\begin{equation} \label{def:cum}
\ka(\pmb{\lambda}/\pmb{\mu}) := (q-1)^{-(\#\pmb{\lambda}/\pmb{\mu} - 1)} \sum\limits_{\mathcal{B}\in\Part(n)} \prod\limits_{B\in\mathcal{B}} (-1)^{|B|-1}(|B|-1)!\LLT(B),
\end{equation}
where $\#\pmb{\lambda}/\pmb{\mu}$ denotes the number of elements in the sequence $\pmb{\lambda}/\pmb{\mu}$, $\Part(n)$ denotes the set of all set partitions of $[n] = \{1,\dots,n\}$, and 
\[\LLT(B):= \LLT(\lambda_{i_1}/\mu_{i_1},\dots,\lambda_{i_m}/\mu_{i_m})\ \  \text{for}\ B=\{i_1 <\cdots < i_m\}.\]

Equivalently, using the M\"obius inversion formula for the poset of set partitions, we can define unicellular LLT cumulants recursively by
\begin{equation} \label{def:cumrecursive}
\LLT(\pmb{\lambda}/\pmb{\mu})(q) = \sum\limits_{\mathcal{B}\in\Part(n)} (q-1)^{n - |\mathcal{B}|} \prod\limits_{B\in\mathcal{B}} \ka(B)(q),
\end{equation}
where
\[\ka(B)(q) := k(\lambda_{i_1}/\mu_{i_1},\dots,\lambda_{i_m}/\mu_{i_m}) \ \  \text{for}\ B=\{i_1 <\cdots < i_m\}.\]

We wish to understand $\LLT(\pmb{\lambda}/\pmb{\mu})$ and $\ka(\pmb{\lambda}/\pmb{\mu})$ as generating functions of certain graphs, which we call \emph{LLT graphs}.

\begin{definition}
Let $G$ be a directed graph with three types of edges, visually depicted as $\rightarrow$, $\twoheadrightarrow$, and $\Rightarrow$, and called \emph{edges of type I}, \emph{of type II}, and \emph{double edges}, respectively. Denote the corresponding sets of edges by $E_1(G)$, $E_2(G)$, and $E_d(G)$. 

\emph{The LLT polynomial corresponding to} $G$ is the expression
\begin{equation} \label{def:graphllt}
\LLT(G):=\sum\limits_{f:V(G)\rightarrow\mathbb{N}} \left(\prod\limits_{(u,v)\in E(G)} \varphi_f(u,v)\right)\cdot \left(\prod\limits_{v\in V(G)} x_{f(v)}\right),
\end{equation}
with
\begin{equation} \label{def:graphcol}
\varphi_f(u,v)=\begin{cases} [f(u)>f(v)] & \text{for } (u,v)\in E_1(G); \\
[f(u)\ge f(v)] & \text{for } (u,v)\in E_2(G); \\
q[f(u)>f(v)]+[f(u)\le f(v)] & \text{for } (u,v)\in E_d(G), \end{cases}
\end{equation}
where $[A]$ is the characteristic function of condition $A$, i.e., is equal to $1$ if $A$ is true and $0$ otherwise.
\end{definition}

\begin{figure}
\centering
\begin{tikzpicture}[scale=.4, every node/.style={scale=0.8}]
\draw[thick, red] (0,0) -- (-1,0) -- (-1,-1) -- (2,-1) -- (2,-2) -- (0,-2) -- (0,0) -- (1,0) -- (1,-2);
\draw[thick, blue] (3,-1) -- (3,0) -- (4,0) -- (4,-2) -- (3,-2) -- (3,-1) -- (4,-1);
\node at (-1,-3) {};
\draw[thick] (2.5,-1.8) arc (0:-45:.5);
\draw[thick] (-1,.5) arc (135:225:2);
\draw[thick] (4,.5) arc (45:-45:2);

\node at (-.5,-.5) {$a$};
\node at (.5,-.5) {$b$};
\node at (.5,-1.5) {$c$};
\node at (1.5,-1.5) {$d$};
\node at (3.5,-.5) {$x$};
\node at (3.5,-1.5) {$y$};
\end{tikzpicture}
\begin{tikzpicture}[scale=1]
\draw[thick, |->] (0,0) -- (1,0);
\node at (-.5,-1) {};
\end{tikzpicture}
\hspace{.5cm}
\begin{tikzpicture}[scale=0.6, every node/.style={scale=0.8}]
\draw[thick, red] (0,0) -- (.5,-.5);
\draw[thick, red, <<-] (.5,-.5) -- (1,-1);
\draw[thick, red, ->] (1,-1) -- (.5,-1.5);
\draw[thick, red] (.5,-1.5) -- (0,-2);
\draw[thick, red] (0,-2) -- (.5,-2.5);
\draw[thick, red, <<-] (.5,-2.5) -- (1,-3);

\draw[thick, blue, ->] (4,-.5) -- (4,-1.5);
\draw[thick, blue] (4,-1.5) -- (4,-2.5);

\draw[thick, double, ->] (0,0) -- (2,-.25);
\draw[thick, double] (2,-.25) -- (4,-.5);
\draw[thick, double, ->] (4,-.5) -- (2.5,-.75);
\draw[thick, double] (2.5,-.75) -- (1,-1);
\draw[thick, double, ->] (1,-1) -- (2.5,-1.75);
\draw[thick, double] (2.5,-1.75) -- (4,-2.5);
\draw[thick, double, ->] (4,-2.5) -- (2,-2.25);
\draw[thick, double] (2,-2.25) -- (0,-2);

\node[below left] at (0,0) {$a$};
\node[above] at (1.1,-1) {$b$};
\node[left] at (0,-2) {$c$};
\node[above right] at (1,-3) {$d$};
\node[right] at (4,-.5) {$x$};
\node[right] at (4,-2.5) {$y$};
\end{tikzpicture}
\caption{The LLT graph corresponding to ((3,2)/(1), (1,1)).}
\label{fig:lltgraph}
\end{figure}
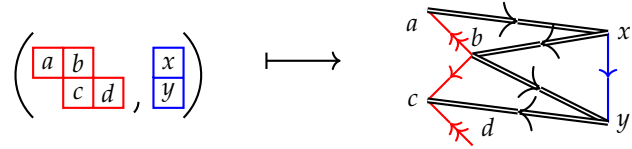

The algebra generated by all linear combinations of LLT graphs over $\mathbb{Z}[q]$ is isomorphic to the algebra of quasi-symmetric functions over $\mathbb{Z}[q]$. Furthermore, there is an obvious way to associate an LLT graph $G_{\pmb{\lambda}/\pmb{\mu}}$ to a sequence of skew shapes $\pmb{\lambda}/\pmb{\mu}$ such that $\LLT(G_{\pmb{\lambda}/\pmb{\mu}}) = LLT(\pmb{\lambda}/\pmb{\mu})$. To be precise, vertices correspond to cells, edges of type I go from a cell (i.e., the vertex associated with that cell) to the cell directly below it, edges of type II go from a cell to that directly to its left, and double edges connect cells that correspond to inversions (see \cref{fig:lltgraph}).

\subsection{Graph relations}

Translating the problem of LLT polynomials into a larger algebra (i.e., from symmetric functions to quasi-symmetric functions) offers much more flexibility. Below, we introduce examples of this in the form of relations on LLT graphs.

\begin{lem} \label{lem:arrowrel}
	\begin{enumerate}[label=(\Alph*)]
		\item $\LLT(\cdot \hspace{.15cm} \cdot)=\LLT(\rightarrow)+\LLT(\twoheadleftarrow)$,
		\item $\LLT(\Rightarrow)=q\LLT(\rightarrow)+\LLT(\twoheadleftarrow)$,
		\item $\LLT(\Rightarrow)=(q-1)\LLT(\rightarrow)+\LLT(\cdot \hspace{.15cm} \cdot)$,
		\item 
		$\LLT\left(\vcenter{\hbox{
				\begin{tikzpicture}[scale=0.7]
				\draw[thick,-{stealth}{stealth}] (-.5,0) -- (-.125,.66);
				\draw[thick] (-.125,.66) -- (0,.88);
				\draw[thick,-{stealth}{stealth}] (0,.88) -- (.375,.22);
				\draw[thick] (.375,.22) -- (.5,0);
				\end{tikzpicture}}}\right)
		=
		\LLT\left(\vcenter{\hbox{
				\begin{tikzpicture}[scale=0.7]
				\draw[thick,-{stealth}{stealth}] (-.5,0) -- (-.125,.66);
				\draw[thick] (-.125,.66) -- (0,.88);
				\draw[thick,-{stealth}{stealth}] (0,.88) -- (.375,.22);
				\draw[thick] (.375,.22) -- (.5,0);
				\draw[thick,-{stealth}{stealth}] (-.5,0) -- (.25,0);
				\draw[thick] (.25,0) -- (.5,0);
				\end{tikzpicture}}}\right)$,
		\item 
		$\LLT\left(\vcenter{\hbox{
				\begin{tikzpicture}[scale=0.7]
				\draw[thick,-{stealth}] (-.5,0) -- (-.125,.66);
				\draw[thick] (-.125,.66) -- (0,.88);
				\draw[thick,-{stealth}] (0,.88) -- (.375,.22);
				\draw[thick] (.375,.22) -- (.5,0);
				\end{tikzpicture}}}\right)
		=
		\LLT\left(\vcenter{\hbox{
				\begin{tikzpicture}[scale=0.7]
				\draw[thick,-{stealth}] (-.5,0) -- (-.125,.66);
				\draw[thick] (-.125,.66) -- (0,.88);
				\draw[thick,-{stealth}] (0,.88) -- (.375,.22);
				\draw[thick] (.375,.22) -- (.5,0);
				\draw[thick,-{stealth}] (-.5,0) -- (.25,0);
				\draw[thick] (.25,0) -- (.5,0);
				\end{tikzpicture}}}\right)$,
		\item If $G$ contains a directed cycle consisting of edges of type I and type II with at least one edge of type I, then $\LLT(G)=0$.
	\end{enumerate}
\end{lem}

Note that due to the multiplicative character of \eqref{def:graphllt}, we can apply any of the above relations locally as long as the rest of the graph remains unchanged.

\subsection{Diagrams corresponding to spanning trees} \label{sec:cumunicell}

From now on, we focus on the case from \cref{thm:singcelldecomp}, i.e., when $\pmb{\lambda}/\pmb{\mu} = ((1),...,(1))$ is a sequence of $n$ unicellular unilevelled diagrams. This corresponds to the full Dyck path $D_n$ in \eqref{def:lltdyck} or a complete graph $K_n$ in \eqref{def:graphllt} with $E(K_n) = E_d(K_n)$.

Let $T$ be a tree on $[n]$. We denote by $\nu(T) = (\nu_1,...,\nu_k)$ the unique sequence of vertical-strip shapes obtained according to the rules described below (for a formal description, see \cite{DolegaKowalski}).

Root $T$ in $1$. Form $\nu_1$ by travelling from $1$ to its child with the smallest label, continuing this way until no longer possible. The number of vertices visited this way is the number of cells in $\nu_1$.

Next, retrace to the previously considered vertices until you reach one with children which have not yet been visited. Choose that with the smallest label, and form $\nu_2$ following the same rules beginning with that vertex.

Repeat the reasoning until you consider all vertices of $T$. Lastly, adjust elements of the sequence so that the cells corresponding to vertices at the same distance from $1$ have the same content in $\nu(T)$ (see the first arrow in \cref{fig:shapeoftree}).

\section{Schr\"{o}der paths and parking functions} \label{sec:schroderpaths}

\textit{A Schr\"{o}der path} of length $n$ is a lattice path from $(0,0)$ to $(n,n)$ with steps $n = (0,1)$, $e = (1,0)$, and $d = (1,1)$ (referred to as \textit{north}, \textit{east}, and \textit{diagonal steps}, respectively), which never falls below the main diagonal that joins the ends, and has no $d$ steps on that diagonal. We denote by $(i,j)$ the coordinates of the $1\times 1$ box with lower left vertex in $(i-1,j-1)$.

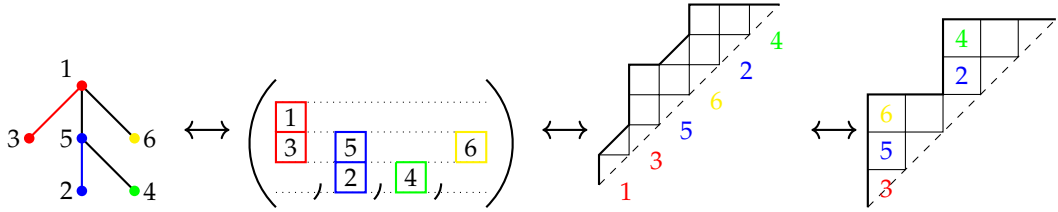
\begin{figure}
	\centering
	\begin{tikzpicture}[scale=.7, every node/.style={scale=0.8}]
	\draw[thick, red] (0,0) -- (-1,-1);
	\draw[thick] (0,0) -- (0,-1);
	\draw[thick] (0,0) -- (1,-1);
	\draw[thick, blue] (0,-1) -- (0,-2);
	\draw[thick] (0,-1) -- (1,-2);
	\fill[red] (0,0) circle (0.1);
	\fill[blue] (0,-1) circle (0.1);
	\fill[blue] (0,-2) circle (0.1);
	\fill[red] (-1,-1) circle (0.1);
	\fill[yellow] (1,-1) circle (0.1);
	\fill[green] (1,-2) circle (0.1);
	\node[above left] at (0,0) {$1$};
	\node[left] at (0,-1) {$5$};
	\node[left] at (0,-2) {$2$};
	\node[left] at (-1,-1) {$3$};
	\node[right] at (1,-1) {$6$};
	\node[right] at (1,-2) {$4$};
	\end{tikzpicture}
	\begin{tikzpicture}[scale=1.2]
	\draw[thick, <->] (0,0) -- (.5,0);
	\node at (0,-.75) {};
	\end{tikzpicture}
	\begin{tikzpicture}[scale=.6, every node/.style={scale=0.8}]
	\draw[dotted] (-.66,0) -- (4,0);
	\draw[dotted] (-.66,-.66) -- (4,-.66);
	\draw[dotted] (-.66,-1.33) -- (4,-1.33);
	\draw[dotted] (-.66,-2) -- (4,-2);
	\draw[thick, red] (-.66,-.66) -- (-.66,0) -- (0,0) -- (0,-1.33) -- (-.66,-1.33) -- (-.66,-.66) -- (0,-.66);
	\draw[thick] (.33,-1.8) arc (0:-45:.5);
	\draw[thick, blue] (.66,-1.33) -- (.66,-.66) -- (1.33,-.66) -- (1.33,-2) -- (.66,-2) -- (.66,-1.33) -- (1.33,-1.33);
	\draw[thick] (1.66,-1.8) arc (0:-45:.5);
	\draw[thick, green] (2,-1.33) -- (2.66,-1.33) -- (2.66,-2) -- (2,-2) -- (2,-1.33);
	\draw[thick] (3,-1.8) arc (0:-45:.5);
	\draw[thick, yellow] (3.33,-1.33) -- (4,-1.33) -- (4,-.66) -- (3.33,-.66) -- (3.33,-1.33);
	\draw[thick] (-.66,.5) arc (135:225:2);
	\draw[thick] (4,.5) arc (45:-45:2);
	\node at (-.33,-.33) {$1$};
	\node at (-.33,-1) {$3$};
	\node at (1,-1) {$5$};
	\node at (1,-1.66) {$2$};
	\node at (2.33,-1.66) {$4$};
	\node at (3.66,-1) {$6$};
	\end{tikzpicture}
	\begin{tikzpicture}[scale=1.2]
	\draw[thick, <->] (0,0) -- (.5,0);
	\node at (0,-.75) {};
	\end{tikzpicture}
	\begin{tikzpicture}[scale=.8, every node/.style={scale=0.8}]
	\draw[dashed] (0,0) -- (3,3);
	\draw (0,.5) -- (.5,.5) -- (.5,1) -- (1,1) -- (1,1.5) -- (1.5,1.5) -- (1.5,2) -- (2,2) -- (2,2.5) -- (2.5,2.5) -- (2.5,3);
	\draw (.5,1.5) -- (1,1.5) -- (1,2) -- (1.5,2) -- (1.5,2.5) -- (2,2.5) -- (2,3);
	\draw[thick] (0,0) -- (0,.5) -- (.5,1) -- (.5,2) -- (1,2) -- (1.5,2.5) -- (1.5,3) -- (3,3);
	\node[red, below right] at (.2,.2) {$1$};
	\node[red, below right] at (.7,.7) {$3$};
	\node[blue, below right] at (1.2,1.2) {$5$};
	\node[yellow, below right] at (1.7,1.7) {$6$};
	\node[blue, below right] at (2.2,2.2) {$2$};
	\node[green, below right] at (2.7,2.7) {$4$};
	\end{tikzpicture}
	\begin{tikzpicture}[scale=1.2]
	\draw[thick, <->] (0,0) -- (.5,0);
	\node at (0,-.75) {};
	\end{tikzpicture}
	\begin{tikzpicture}[scale=1, every node/.style={scale=0.8}]
	\draw[dashed] (.5,.5) -- (3,3);
	\draw (.5,.5) -- (.5,1) -- (1,1) -- (1,1.5) -- (1.5,1.5) -- (1.5,2) -- (2,2) -- (2,2.5) -- (2.5,2.5) -- (2.5,3);
	\draw (.5,1.5) -- (1,1.5) -- (1,2) -- (1.5,2) -- (1.5,2.5) -- (2,2.5) -- (2,3);
	\draw[thick] (.5,.5) -- (.5,2) -- (1.5,2) -- (1.5,3) -- (3,3);
	\node[red] at (.75,.75) {$3$};
	\node[blue] at (.75,1.25) {$5$};
	\node[yellow] at (.75,1.75) {$6$};
	\node[blue] at (1.75,2.25) {$2$};
	\node[green] at (1.75,2.75) {$4$};
	\end{tikzpicture}
	\caption{A spanning tree with its corresponding sequence of shapes, Schr\"{o}der path, and parking function.}
	\label{fig:shapeoftree}
\end{figure}

For such a path $P$, we say that its $i$-th column \textit{is of height} $h(i)=k$ if $(i-1,k)$ lies on $P$ and is followed by either an east or diagonal step, $1\le i\le n$. We define the \emph{jump} of column $i$ to be the value $j(i)=h(i+1)-h(i)$. Lastly, a box $(i,j)$ is called \textit{an outer corner} of $P$ if the point $(i-1,j-1)$ lies on $P$ and is followed by an east step and a north step.

The sequences of vertical-strip shapes are in bijection with Schr\"oder paths. To be precise, for a path $P$ of length $n$, label box $(i,i)$ by the integer $i$, $1\le i\le n$. Then, construct the sequence $\mu(P)=(\mu_1,...,\mu_k)$ of vertical-strip shapes according to the following rules (for precise formulation, see \cite{DolegaKowalski}):
\begin{enumerate}
	\item if $P$ has a diagonal step in the box $(i,j)$, then cell $i$ is directly above cell $j$ in some $\mu_m$;
	\item if the box $(i,j)$, $i<j$, lies under $P$, then the pair of cells $(i,j)$ is an inversion pair in $\mu(P)$.
\end{enumerate}

The bijection given in \cref{sec:cumunicell} and the above correspondence induces an injection from planar rooted trees with $n$ vertices to Schr\"{o}der paths of length $n$: $T \leftrightarrow \nu(T) = \mu(P) \leftrightarrow P$. This injection turns into bijection after reducing the set of Schr\"oder paths to those which satisfy the following conditions (see \cref{fig:shapeoftree}):
\begin{enumerate}
	\item $P$ is connected (i.e., the only points of $P$ on the diagonal are $(0,0)$ and $(n,n)$);
	\item the first steps of $P$ are $n$ and $d$; and
	\item $P$ has no outer corners.
\end{enumerate}

Finally, it turns out that the above bijection induces a bijection between Cayley trees on $n$ vertices and parking functions of length $n-1$. Recall that \emph{a parking function} $f:[n-1]\rightarrow[n-1]$ is a map for which $|f^{-1}([i])|\ge i$ for $1\le i\le n-1$. We can represent such objects geometrically as Dyck paths between $(1,1)$ and $(n,n)$ (i.e., Schr\"oder paths with no diagonal steps) with boxes to the right of north steps labeled with integers $2$ through $n$ so that the numbers increase in columns upwards. We denote by $P_f$ the labeled Schr\"oder path of length $n$ obtained from $f$ by adding a cell $(1,1)$ labeled by $1$ and replacing all steps $en$ by $d$ (see \cref{fig:shapeoftree}). We then define $\mu(f):=\mu(P_f)$ and $\nu^{-1}\circ \mu$ is the desired bijection.

\subsection{Schr\"oder path relations}

Similarly to how we introduced relations between graphs (or rather between the generating functions of the graphs), we can study relations between Schr\"{o}der paths. This approach was used by Alexandersson and Sulzgruber \cite{AlexSulz}, who show the following

\begin{lem} \label{lem:schroderrelations}
Let $P$ be a Schr\"{o}der path.
\begin{enumerate}[label=(\Alph*)]
	\item If $P=SneT$ for some paths $S$ and $T$, then
	\[\LLT(SneT) = (q-1)\LLT(SdT) + \LLT(SenT).\]
	\item If $P = SndReeT$ for some paths $S$, $T$, and $R$ with $S$ ending in $(i,j)$ and $SndR$ ending in $(j,k)$, then
	\[\LLT(SndReeT) = q\LLT(SdnReeT).\]
\end{enumerate}
\end{lem}

The visual representation of the above relations in the spirit of \cref{lem:arrowrel} takes the form
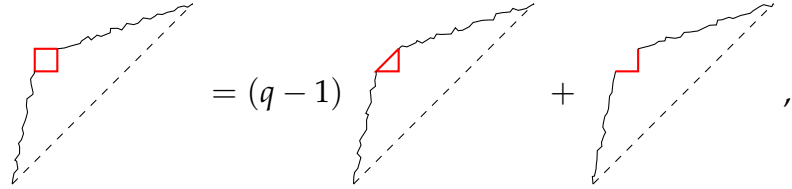
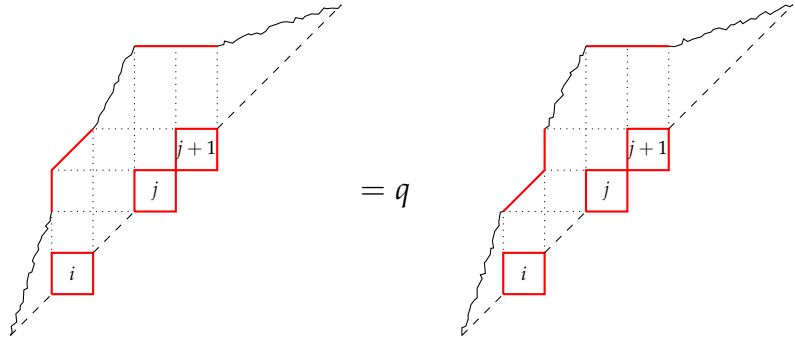
\begin{enumerate}[label=(\Alph*)]
	\item
\begin{equation*}
\vcenter{
	\hbox{
		\begin{tikzpicture}[scale=0.6]
		\draw[thick, red] (1,3) -- (.5,3) -- (.5,2.5) -- (1,2.5) -- (1,3);
		\draw[dashed] (0,0) -- (4,4);
		\path[draw,clip,decoration={random steps,segment length=2pt,amplitude=1pt}] decorate {(1,3) -- (4,4)} decorate {(0,0) -- (0.5,2.5)};
		\end{tikzpicture}
	}
}
=
(q-1)\vcenter{
	\hbox{
		\begin{tikzpicture}[scale=0.6]
		\draw[thick, red] (1,3) -- (.5,2.5) -- (1,2.5) -- (1,3);
		\draw[dashed] (0,0) -- (4,4);
		\path[draw,clip,decoration={random steps,segment length=2pt,amplitude=1pt}] decorate {(1,3) -- (4,4)} decorate {(0,0) -- (0.5,2.5)};
		\end{tikzpicture}
	}
}
+
\vcenter{
	\hbox{
		\begin{tikzpicture}[scale=0.6]
		\draw[thick, red] (1,3) -- (1,2.5) -- (.5,2.5);
		\draw[dashed] (0,0) -- (4,4);
		\path[draw,clip,decoration={random steps,segment length=2pt,amplitude=1pt}] decorate {(1,3) -- (4,4)} decorate {(0,0) -- (0.5,2.5)};
		\end{tikzpicture}
	}
},
\end{equation*}

\item
\begin{equation*}
\vcenter{
	\hbox{
		\begin{tikzpicture}[scale=1.1, every node/.style={scale=0.6}]
		\draw[dotted] (1,2.5) -- (2.5,2.5) -- (2.5,4);
		\draw[dotted] (1,1) -- (1,2) -- (2,2) -- (2,4);
		\draw[dotted] (1.5,1.5) -- (1.5,3) -- (3,3) -- (3,4);
		\draw[thick, red] (1,2) -- (1,2.5) -- (1.5,3);
		\draw[thick, red] (2,4) -- (3,4);
		\draw[thick, red] (1,1) -- (1,1.5) -- (1.5,1.5) -- (1.5,1) -- (1,1);
		\draw[thick, red] (2,2) -- (2,2.5) -- (3,2.5) -- (3,3) -- (2.5,3) -- (2.5,2) -- (2,2);
		\draw[dashed] (.5,.5) -- (1,1);
		\draw[dashed] (1.5,1.5) -- (2,2);
		\draw[dashed] (3,3) -- (4.5,4.5);
		\node at (1.25,1.25) {$i$};
		\node at (2.25,2.25) {$j$};
		\node at (2.75,2.75) {$j+1$};
		\path[draw,clip,decoration={random steps,segment length=2pt,amplitude=1pt}] decorate {(3,4) -- (4.5,4.5)} decorate {(.5,.5) -- (1,2)} decorate {(1.5,3) -- (2,4)};
		\end{tikzpicture}
	}
}
=
q\vcenter{
	\hbox{
		\begin{tikzpicture}[scale=1.1, every node/.style={scale=0.6}]
		\draw[dotted] (1.5,2.5) -- (2.5,2.5) -- (2.5,4);
		\draw[dotted] (1,1) -- (1,2) -- (2,2) -- (2,4);
		\draw[dotted] (1.5,1.5) -- (1.5,3) -- (3,3) -- (3,4);
		\draw[thick, red] (1,2) -- (1.5,2.5) -- (1.5,3);
		\draw[thick, red] (2,4) -- (3,4);
		\draw[thick, red] (1,1) -- (1,1.5) -- (1.5,1.5) -- (1.5,1) -- (1,1);
		\draw[thick, red] (2,2) -- (2,2.5) -- (3,2.5) -- (3,3) -- (2.5,3) -- (2.5,2) -- (2,2);
		\draw[dashed] (.5,.5) -- (1,1);
		\draw[dashed] (1.5,1.5) -- (2,2);
		\draw[dashed] (3,3) -- (4.5,4.5);
		\node at (1.25,1.25) {$i$};
		\node at (2.25,2.25) {$j$};
		\node at (2.75,2.75) {$j+1$};
		\path[draw,clip,decoration={random steps,segment length=2pt,amplitude=1pt}] decorate {(3,4) -- (4.5,4.5)} decorate {(.5,.5) -- (1,2)} decorate {(1.5,3) -- (2,4)};
		\end{tikzpicture}
	}
}.
\end{equation*}
\end{enumerate}
Additionally, the second identity locally (i.e., limited to the three vertices corresponding to the marked boxes) translates to the following identity on LLT graphs:
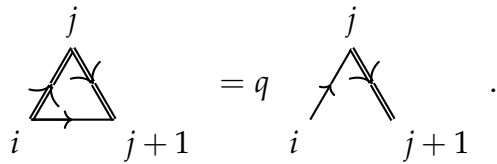
\begin{equation*}
\vcenter{
	\hbox{
		\begin{tikzpicture}[scale=1.1]
		\draw[thick,double, ->] (-.5,0) -- (-.25,.43);
		\draw[thick,double] (-.25,.43) -- (0,.86);
		\draw[thick,double, ->] (0,.86) -- (.25,.43);
		\draw[thick,double] (.25,.43) -- (.5,0);
		\draw[thick, ->] (-0.5,0) -- (0,0);
		\draw[thick] (-0.5,0) -- (0.5,0);
		\node[below left] at (-0.5,0) {$i$};
		\node[above] at (0,0.86) {$j$};
		\node[below right] at (0.5,0) {$j+1$};
		\end{tikzpicture}
	}
}
=q
\vcenter{
	\hbox{
		\begin{tikzpicture}[scale=1.1]
		\draw[thick, ->] (-.5,0) -- (-.25,.43);
		\draw[thick] (-.25,.43) -- (0,.86);
		\draw[thick,double, ->] (0,.86) -- (.25,.43);
		\draw[thick,double] (.25,.43) -- (.5,0);
		\node[below left] at (-0.5,0) {$i$};
		\node[above] at (0,0.86) {$j$};
		\node[below right] at (0.5,0) {$j+1$};
		\end{tikzpicture}
	}
}.
\end{equation*}

We present a different proof of \cref{lem:schroderrelations} from the one in \cite{AlexSulz}. It is an easy manipulation of relations presented in \cref{lem:arrowrel}.

\begin{proof}
\begin{enumerate}[label=(\Alph*)]
	\item The identity is a straightforward translation of property $(3)$ from \cref{lem:arrowrel}
	
	\item From the shape of $P$ we deduce that the only relations that change after transformation are between cells $i$, $j$, and $j+1$. In other words, it is enough to look at the objects locally.
	
	If we use graph relations from \cref{lem:arrowrel}, we get (for simplicity of notation, we drop the vertex labels)
	\begin{align*}
		\vcenter{\hbox{
				\begin{tikzpicture}[scale=1.1]
				\draw[thick,double, ->] (-.5,0) -- (-.25,.43);
				\draw[thick,double] (-.25,.43) -- (0,.86);
				\draw[thick,double, ->] (0,.86) -- (.25,.43);
				\draw[thick,double] (.25,.43) -- (.5,0);
				\draw[thick, ->] (-0.5,0) -- (0,0);
				\draw[thick] (-0.5,0) -- (0.5,0);
				\end{tikzpicture}}}
		&\hspace{.4cm}\stackrel{(C)}{=}\hspace{.4cm}
		(q-1)\vcenter{\hbox{
				\begin{tikzpicture}[scale=1.1]
				\draw[thick, ->] (-.5,0) -- (-.25,.43);
				\draw[thick] (-.25,.43) -- (0,.86);
				\draw[thick,double, ->] (0,.86) -- (.25,.43);
				\draw[thick,double] (.25,.43) -- (.5,0);
				\draw[thick, ->] (-0.5,0) -- (0,0);
				\draw[thick] (-0.5,0) -- (0.5,0);
				\end{tikzpicture}}}
		+
		\vcenter{\hbox{
				\begin{tikzpicture}[scale=1.1]
				\draw[thick,double, ->] (0,.86) -- (.25,.43);
				\draw[thick,double] (.25,.43) -- (.5,0);
				\draw[thick, ->] (-0.5,0) -- (0,0);
				\draw[thick] (-0.5,0) -- (0.5,0);
				\end{tikzpicture}}}
		= \\
		&\hspace{.4cm}\stackrel{(A)}{=}\hspace{.4cm}
		(q-1)\vcenter{\hbox{
				\begin{tikzpicture}[scale=1.1]
				\draw[thick, ->] (-.5,0) -- (-.25,.43);
				\draw[thick] (-.25,.43) -- (0,.86);
				\draw[thick,double, ->] (0,.86) -- (.25,.43);
				\draw[thick,double] (.25,.43) -- (.5,0);
				\end{tikzpicture}}}
		-
		(q-1)\vcenter{\hbox{
				\begin{tikzpicture}[scale=1.1]
				\draw[thick, ->] (-.5,0) -- (-.25,.43);
				\draw[thick] (-.25,.43) -- (0,.86);
				\draw[thick,double, ->] (0,.86) -- (.25,.43);
				\draw[thick,double] (.25,.43) -- (.5,0);
				\draw[thick] (-0.5,0) -- (-.1,0);
				\draw[thick, <<-] (-0.1,0) -- (0.5,0);
				\end{tikzpicture}}}
		+
		\vcenter{\hbox{
				\begin{tikzpicture}[scale=1.1]
				\draw[thick,double, ->] (0,.86) -- (.25,.43);
				\draw[thick,double] (.25,.43) -- (.5,0);
				\filldraw[black] (-0.5,0) circle (1pt);
				\end{tikzpicture}}}
		-
		\vcenter{\hbox{
				\begin{tikzpicture}[scale=1.1]
				\draw[thick,double, ->] (0,.86) -- (.25,.43);
				\draw[thick,double] (.25,.43) -- (.5,0);
				\draw[thick] (-0.5,0) -- (-.1,0);
				\draw[thick, <<-] (-0.1,0) -- (0.5,0);
				\end{tikzpicture}}}
		= \hspace{1cm} \\
		&\stackrel{(C)+(A)}{=}
		(q-1)\vcenter{\hbox{
				\begin{tikzpicture}[scale=1.1]
				\draw[thick, ->] (-.5,0) -- (-.25,.43);
				\draw[thick] (-.25,.43) -- (0,.86);
				\draw[thick,double, ->] (0,.86) -- (.25,.43);
				\draw[thick,double] (.25,.43) -- (.5,0);
				\end{tikzpicture}}}
		-
		(q-1)\vcenter{\hbox{
				\begin{tikzpicture}[scale=1.1]
				\draw[thick, ->] (-.5,0) -- (-.25,.43);
				\draw[thick] (-.25,.43) -- (0,.86);
				\draw[thick, ->] (0,.86) -- (.25,.43);
				\draw[thick] (.25,.43) -- (.5,0);
				\draw[thick] (-0.5,0) -- (-.1,0);
				\draw[thick, <<-] (-0.1,0) -- (0.5,0);
				\end{tikzpicture}}}
		-
		(q-1)\vcenter{\hbox{
				\begin{tikzpicture}[scale=1.1]
				\draw[thick, ->] (-.5,0) -- (-.25,.43);
				\draw[thick] (-.25,.43) -- (0,.86);
				\draw[thick] (-0.5,0) -- (-.1,0);
				\draw[thick, <<-] (-0.1,0) -- (0.5,0);
				\end{tikzpicture}}}
		+ \\
		&\hspace{.4cm}+\hspace{.4cm}
		\vcenter{\hbox{
				\begin{tikzpicture}[scale=1.1]
				\draw[thick, ->] (-.5,0) -- (-.25,.43);
				\draw[thick] (-.25,.43) -- (0,.86);
				\draw[thick,double, ->] (0,.86) -- (.25,.43);
				\draw[thick,double] (.25,.43) -- (.5,0);
				\end{tikzpicture}}}
		+
		\vcenter{\hbox{
				\begin{tikzpicture}[scale=1.1]
				\draw[thick] (-.5,0) -- (-.3,.336);
				\draw[thick, <<-] (-.3,.336) -- (0,.86);
				\draw[thick,double, ->] (0,.86) -- (.25,.43);
				\draw[thick,double] (.25,.43) -- (.5,0);
				\end{tikzpicture}}}
		-
		\vcenter{\hbox{
				\begin{tikzpicture}[scale=1.1]
				\draw[thick,double, ->] (0,.86) -- (.25,.43);
				\draw[thick,double] (.25,.43) -- (.5,0);
				\draw[thick] (-0.5,0) -- (-.1,0);
				\draw[thick, <<-] (-0.1,0) -- (0.5,0);
				\end{tikzpicture}}}
		= \\
		&\stackrel{(F)+(C)}{=}
		q\vcenter{\hbox{
				\begin{tikzpicture}[scale=1.1]
				\draw[thick, ->] (-.5,0) -- (-.25,.43);
				\draw[thick] (-.25,.43) -- (0,.86);
				\draw[thick,double, ->] (0,.86) -- (.25,.43);
				\draw[thick,double] (.25,.43) -- (.5,0);
				\end{tikzpicture}}}
		-
		(q-1)\vcenter{\hbox{
				\begin{tikzpicture}[scale=1.1]
				\draw[thick, ->] (-.5,0) -- (-.25,.43);
				\draw[thick] (-.25,.43) -- (0,.86);
				\draw[thick] (-0.5,0) -- (-.1,0);
				\draw[thick, <<-] (-0.1,0) -- (0.5,0);
				\end{tikzpicture}}}
		+ \\
		&\hspace{.4cm}+\hspace{.4cm}
		(q-1)\vcenter{\hbox{
				\begin{tikzpicture}[scale=1.1]
				\draw[thick] (-.5,0) -- (-.3,.336);
				\draw[thick, <<-] (-.3,.336) -- (0,.86);
				\draw[thick, ->] (0,.86) -- (.25,.43);
				\draw[thick] (.25,.43) -- (.5,0);
				\end{tikzpicture}}}
		+
		\vcenter{\hbox{
				\begin{tikzpicture}[scale=1.1]
				\draw[thick] (-.5,0) -- (-.3,.336);
				\draw[thick, <<-] (-.3,.336) -- (0,.86);
				\filldraw[black] (0.5,0) circle (1pt);
				\end{tikzpicture}}}
		-
		(q-1)\vcenter{\hbox{
				\begin{tikzpicture}[scale=1.1]
				\draw[thick, ->] (0,.86) -- (.25,.43);
				\draw[thick] (.25,.43) -- (.5,0);
				\draw[thick] (-0.5,0) -- (-.1,0);
				\draw[thick, <<-] (-0.1,0) -- (0.5,0);
				\end{tikzpicture}}}
		-
		\vcenter{\hbox{
				\begin{tikzpicture}[scale=1.1]
				\draw[thick] (-0.5,0) -- (-.1,0);
				\draw[thick, <<-] (-0.1,0) -- (0.5,0);
				\filldraw[black] (0,0.86) circle (1pt);
				\end{tikzpicture}}}
		= \\
		&\hspace{.4cm}\stackrel{(A)}{=}\hspace{.4cm}
		q\vcenter{\hbox{
				\begin{tikzpicture}[scale=1.1]
				\draw[thick, ->] (-.5,0) -- (-.25,.43);
				\draw[thick] (-.25,.43) -- (0,.86);
				\draw[thick,double, ->] (0,.86) -- (.25,.43);
				\draw[thick,double] (.25,.43) -- (.5,0);
				\end{tikzpicture}}}
		-
		(q-1)\vcenter{\hbox{
				\begin{tikzpicture}[scale=1.1]
				\draw[thick, ->] (-.5,0) -- (-.25,.43);
				\draw[thick] (-.25,.43) -- (0,.86);
				\draw[thick] (-0.5,0) -- (-.1,0);
				\draw[thick, <<-] (-0.1,0) -- (0.5,0);
				\end{tikzpicture}}}
		+
		(q-1)\vcenter{\hbox{
				\begin{tikzpicture}[scale=1.1]
				\draw[thick] (-.5,0) -- (-.3,.336);
				\draw[thick, <<-] (-.3,.336) -- (0,.86);
				\draw[thick, ->] (0,.86) -- (.25,.43);
				\draw[thick] (.25,.43) -- (.5,0);
				\draw[thick, ->] (-0.5,0) -- (0,0);
				\draw[thick] (-0.5,0) -- (0.5,0);
				\end{tikzpicture}}}
		+
		(q-1)\vcenter{\hbox{
				\begin{tikzpicture}[scale=1.1]
				\draw[thick] (-.5,0) -- (-.3,.336);
				\draw[thick, <<-] (-.3,.336) -- (0,.86);
				\draw[thick, ->] (0,.86) -- (.25,.43);
				\draw[thick] (.25,.43) -- (.5,0);
				\draw[thick] (-0.5,0) -- (-.1,0);
				\draw[thick, <<-] (-0.1,0) -- (0.5,0);
				\end{tikzpicture}}}
		+ \\
		&\hspace{.4cm}+\hspace{.4cm}
		\vcenter{\hbox{
				\begin{tikzpicture}[scale=1.1]
				\draw[thick] (-.5,0) -- (-.3,.336);
				\draw[thick, <<-] (-.3,.336) -- (0,.86);
				\filldraw[black] (0.5,0) circle (1pt);
				\end{tikzpicture}}}
		-
		(q-1)\vcenter{\hbox{
				\begin{tikzpicture}[scale=1.1]
				\draw[thick, ->] (0,.86) -- (.25,.43);
				\draw[thick] (.25,.43) -- (.5,0);
				\draw[thick] (-0.5,0) -- (-.1,0);
				\draw[thick, <<-] (-0.1,0) -- (0.5,0);
				\end{tikzpicture}}}
		-
		\vcenter{\hbox{
				\begin{tikzpicture}[scale=1.1]
				\draw[thick] (-0.5,0) -- (-.1,0);
				\draw[thick, <<-] (-0.1,0) -- (0.5,0);
				\filldraw[black] (0,0.86) circle (1pt);
				\end{tikzpicture}}}
		= \\
		&\hspace{.4cm}=\hspace{.4cm}
		q\vcenter{\hbox{
				\begin{tikzpicture}[scale=1.1]
				\draw[thick, ->] (-.5,0) -- (-.25,.43);
				\draw[thick] (-.25,.43) -- (0,.86);
				\draw[thick,double, ->] (0,.86) -- (.25,.43);
				\draw[thick,double] (.25,.43) -- (.5,0);
				\end{tikzpicture}}}.
	\end{align*}
	The last equality follows from
	\begin{equation*}
	\vcenter{\hbox{
			\begin{tikzpicture}[scale=1.1]
			\draw[thick] (-.5,0) -- (-.3,.336);
			\draw[thick, <<-] (-.3,.336) -- (0,.86);
			\draw[thick, ->] (0,.86) -- (.25,.43);
			\draw[thick] (.25,.43) -- (.5,0);
			\draw[thick] (-0.5,0) -- (-.1,0);
			\draw[thick, <<-] (-0.1,0) -- (0.5,0);
			\end{tikzpicture}}}
	-
	\vcenter{\hbox{
			\begin{tikzpicture}[scale=1.1]
			\draw[thick, ->] (0,.86) -- (.25,.43);
			\draw[thick] (.25,.43) -- (.5,0);
			\draw[thick] (-0.5,0) -- (-.1,0);
			\draw[thick, <<-] (-0.1,0) -- (0.5,0);
			\end{tikzpicture}}}
	=
	0
	\end{equation*}
	(which is a result of \cref{lem:arrowrel}) and
	\begin{equation*}
	\vcenter{\hbox{
			\begin{tikzpicture}[scale=1.1]
			\draw[thick] (-.5,0) -- (-.3,.336);
			\draw[thick, <<-] (-.3,.336) -- (0,.86);
			\filldraw[black] (0.5,0) circle (1pt);
			\end{tikzpicture}}}
	-
	\vcenter{\hbox{
			\begin{tikzpicture}[scale=1.1]
			\draw[thick] (-0.5,0) -- (-.1,0);
			\draw[thick, <<-] (-0.1,0) -- (0.5,0);
			\filldraw[black] (0,0.86) circle (1pt);
			\end{tikzpicture}}}
	=
	0,
	\hspace{2cm}
	\vcenter{\hbox{
			\begin{tikzpicture}[scale=1.1]
			\draw[thick, ->] (-.5,0) -- (-.25,.43);
			\draw[thick] (-.25,.43) -- (0,.86);
			\draw[thick] (-0.5,0) -- (-.1,0);
			\draw[thick, <<-] (-0.1,0) -- (0.5,0);
			\end{tikzpicture}}}
	-
	\vcenter{\hbox{
			\begin{tikzpicture}[scale=1.1]
			\draw[thick] (-.5,0) -- (-.3,.336);
			\draw[thick, <<-] (-.3,.336) -- (0,.86);
			\draw[thick, ->] (0,.86) -- (.25,.43);
			\draw[thick] (.25,.43) -- (.5,0);
			\draw[thick, ->] (-0.5,0) -- (0,0);
			\draw[thick] (-0.5,0) -- (0.5,0);
			\end{tikzpicture}}}
	=
	0
	\end{equation*}
	which results from the following symmetric property of our graph: for every $v\in V\setminus\{i,j,j+1\}$, the pair $(j,v)$ is an edge of type $P\in\{\rightarrow,\twoheadrightarrow,\Rightarrow\}$ if and only if $(j+1,v)$ is an edge of the same type $P$. Therefore, we can exchange vertices $j$ and $j+1$.
\end{enumerate}
\end{proof}

\section{Proof of \cref{thm:singcelldecomp}} \label{sec:proof}

\begin{proof}[Sketch of the proof]
We are going to prove a formula equivalent to \eqref{thm:formulawithtrees}:
\begin{equation} \label{thm:formulawithforests}
\LLT_{D_n}(q+1) = \sum\limits_F q^{n - \#F}\LLT(\nu(F))(q+1),
\end{equation}
where the sum runs over all forests on $[n]$, $\#F$ denotes the number of connected components of $F$, and $\LLT(\nu(F))$ is the product of $\LLT(\nu(T_i))$ with $T_i$ the connected components of $F$, each rooted in its vertex of the smallest label. The equivalence of the two formulas results from \eqref{def:cumrecursive}.

The idea behind the proof is to take the full Schr\"oder path $D_n$ and repeatedly apply \cref{lem:schroderrelations} to decompose it into a sum of Schr\"oder paths corresponding to forests on $[n]$.

We begin by reducing the height of the first column in $D_n$. First, we apply relation $(A)$ from \cref{lem:schroderrelations} to the cell $(1,n)$, which gives
\begin{equation} \label{thm:reduce1stcolumn}
\vcenter{
	\hbox{
		\begin{tikzpicture}[scale=.5]
		\draw[thick] (0,0) -- (0,3) -- (3,3);
		\draw (0,2.5) -- (0.5,2.5) -- (0.5,3);
		\draw[dashed] (0,0) -- (3,3);
		\end{tikzpicture}
	}
}
=
q\vcenter{
	\hbox{
		\begin{tikzpicture}[scale=.5]
		\draw[thick] (0,0) -- (0,2.5) -- (0.5,3) -- (3,3);
		\draw (0,2.5) -- (0.5,2.5) -- (0.5,3);
		\draw[dashed] (0,0) -- (3,3);
		\end{tikzpicture}
	}
}
+
\vcenter{
	\hbox{
		\begin{tikzpicture}[scale=.5]
		\draw[thick] (0,0) -- (0,2.5) -- (0.5,2.5) -- (0.5,3) -- (3,3);
		\draw[dashed] (0,0) -- (3,3);
		\end{tikzpicture}
	}
}.
\end{equation}

Next, to the first summand, we apply relation $(B)$ from \cref{lem:schroderrelations} to move the diagonal step downwards and repeat until the diagonal step is in the cell $(1,2)$.

To the second summand from \eqref{thm:reduce1stcolumn}, we apply $(A)$ from \cref{lem:schroderrelations} to cell $(1,n-1)$. This will again give us two summands, one with a diagonal step in the first column and one without. To these two, we repeat the reasoning presented.

Eventually, we end up with
\begin{equation} \label{thm:firstcolumndone}
\vcenter{
	\hbox{
		\begin{tikzpicture}[scale=.5]
		\draw[thick] (0,0) -- (0,3) -- (3,3);
		\draw[dashed] (0,0) -- (3,3);
		\end{tikzpicture}
	}
}
=
P(h(1);q)\vcenter{
	\hbox{
		\begin{tikzpicture}[scale=.5]
		\draw[thick] (0,0) -- (0,0.5) -- (0.5,1) -- (0.5,3) -- (3,3);
		\draw (0,0.5) -- (0.5,0.5) -- (0.5,1);
		\draw[dashed] (0,0) -- (3,3);
		\end{tikzpicture}
	}
}
+
\vcenter{
	\hbox{
		\begin{tikzpicture}[scale=.5]
		\draw[thick] (0,0) -- (0,0.5) -- (0.5,0.5) -- (0.5,3) -- (3,3);
		\draw[dashed] (0,0) -- (3,3);
		\end{tikzpicture}
	}
},
\end{equation}
where
\begin{equation*}
P(h(1);q) = q(q+1)^{n-2} + q(q+1)^{n-3} + \cdots + q = (q+1)^{n-1} - 1 = \sum_{i=1}^{n-1} {n-h(1) \choose i}q^i.
\end{equation*}

The next steps of the decomposition differ for every summand. For the second shape in \eqref{thm:firstcolumndone}, we reduce the second column as we did above. For the first, for $1\le k\le n-1$, we decompose the summand with $q^k$ until the second column is of height $k+1$ (equivalently, until $j(1)=k$). In other words, we use the two identities from \cref{lem:schroderrelations} to get
\begin{align*}
{n-h(1) \choose n-h(2)}q^{j(1)}&\vcenter{
	\hbox{
		\begin{tikzpicture}[scale=.5]
		\draw[thick] (0,0) -- (0,0.5) -- (0.5,1) -- (0.5,3) -- (3,3);
		\draw (0,0) -- (0,0.5) -- (0.5,0.5) -- (0.5,1);
		\draw[dashed] (0,0) -- (3,3);
		\end{tikzpicture}
	}
}
= \\
&=
{n-h(1) \choose n-h(2)}q^{j(1)}P(h(2);q)\vcenter{
	\hbox{
		\begin{tikzpicture}[scale=.5]
		\draw[thick] (0,0) -- (0,0.5) -- (0.5,1) -- (0.5,2) -- (1,2.5) -- (1,3) -- (3,3);
		\draw (0.5,2) -- (1,2) -- (1,2.5);
		\draw[dashed] (0,0) -- (3,3);
		\end{tikzpicture}
	}
}
+
{n-h(1) \choose n-h(2)}q^{j(1)} \vcenter{
	\hbox{
		\begin{tikzpicture}[scale=.5]
		\draw[thick] (0,0) -- (0,0.5) -- (0.5,1) -- (0.5,2) -- (1,2) -- (1,3) -- (3,3);
		\draw[dashed] (0,0) -- (3,3);
		\end{tikzpicture}
	}
},
\end{align*}
where, analogously to the previous step,
\begin{equation*}
P(h(2);q) = (q+1)^{n-h(2)} - 1 = \sum_{i=1}^{n-h(2)} {n-h(2) \choose i}q^i.
\end{equation*}

In general, assume that we have decomposed the first $2\le m\le n-1$ columns, i.e., that we have managed to express the left-hand side of \eqref{thm:formulawithforests} as a sum of Schr\"{o}der paths whose first $m$ columns have been reduced.

Take a summand with coefficient $q^K$ of a path $P$ with diagonal steps in columns $1\le i_1,\dots,i_l= m$, $K=j(i_1)+\cdots+j(i_{l-1})+k$, where $1 \leq k \leq n-h(m)$ (if $i_l < m$ then $K=j(i_1)+\cdots+j(i_l)$ and we reduce the $(m+1)$-th column as in the first step). We lower the $(m+1)$-th column to height $K+1$.

\begin{align*}
&\prod_{s=1}^{l-1}{ n-h(i_s) \choose n- h(i_s+1)}\cdot { n-h(i_l) \choose k} q^K\vcenter{
	\hbox{
		\begin{tikzpicture}[scale=.5]
		\draw[thick] (0.5,1.5) -- (1,2) -- (1,4) -- (4,4);
		\draw[dashed] (0,0) -- (4,4);
		\path[draw,clip,decoration={random steps,segment length=2pt,amplitude=1pt}] decorate {(0,0) -- (0.5,1.5)};
		\end{tikzpicture}
	}
}
= \\
		&=\prod_{s=1}^{l}{ n-h(i_s) \choose n- h(i_s+1)} q^KP(h(i_l+1);q)\vcenter{
	\hbox{
		\begin{tikzpicture}[scale=.5]
		\draw[thick] (0.5,1.5) -- (1,2) -- (1,2.5) -- (1.5,3) -- (1.5,4) -- (4,4);
		\draw (1,2.5) -- (1.5,2.5) -- (1.5,3);
		\draw[dashed] (0,0) -- (4,4);
		\path[draw,clip,decoration={random steps,segment length=2pt,amplitude=1pt}] decorate {(0,0) -- (0.5,1.5)};
		\end{tikzpicture}
	}
}
+\prod_{s=1}^{l}{ n-h(i_s) \choose n- h(i_s+1)} q^K\vcenter{
	\hbox{
		\begin{tikzpicture}[scale=.5]
		\draw[thick] (0.5,1.5) -- (1,2) -- (1,2.5) -- (1.5,2.5) -- (1.5,4) -- (4,4);
		\draw[dashed] (0,0) -- (4,4);
		\path[draw,clip,decoration={random steps,segment length=2pt,amplitude=1pt}] decorate {(0,0) -- (0.5,1.5)};
		\end{tikzpicture}
	}
},
\end{align*}
where $K = j(i_1)+\cdots+j(i_l)$ (meaning that $k = j(i_l)$) and 
\begin{equation*}
P(h(i_l+1);q) = (q+1)^{n-h(i_l+1)} - 1 = \sum_{i=1}^{n-h(i_l+1)} {n-h(i_l+1) \choose i}q^i.
\end{equation*}

In the end, the above decomposition gives
\begin{equation} \label{thm:decompositioninpaths}
\LLT_{D_n}(q+1) = \sum_P \varphi(P;q)\LLT(\mu(P))(q+1),
\end{equation}
where the sum runs over all Schr\"oder paths $P$ whose each connected component of length at least $2$ (i.e., a section of the path from one point on the diagonal to another without any others in between) begins with $nd$ and has no outer corners, and
\begin{equation*}
\varphi(P;q) = \prod_{s=1}^{l}{ n-h(i_s) \choose n- h(i_s+1)},
\end{equation*}
where $i_1,\dots,i_l$ are the columns of $P$ containing diagonal steps.

Therefore, each shape appearing in \eqref{thm:decompositioninpaths} corresponds to a forest in the sense of \cref{sec:schroderpaths}. Furthermore, a path $P$ having diagonal steps in columns $i_1,\dots,i_l$ corresponds to the vertex $i_k$ having $j(i_k)$ children in the forest. It is an easy exercise to check that the number of such forests is given by the multinomial coefficient $\varphi(P;q)$, which finishes the proof.
\end{proof}


\printbibliography

\end{document}